\documentclass{article}

\usepackage{amssymb, amsmath, amsthm, enumitem, multirow, cancel, physics, color, mathtools, xcolor, mathrsfs , hyperref,tkz-euclide, tikz, bigints}

\usepackage{mathtools}
\DeclarePairedDelimiter\ceil{\lceil}{\rceil}

\usepackage{stackengine}

\usepackage{arxiv}

\usepackage[utf8]{inputenc} 
\usepackage[T1]{fontenc}    
\usepackage{hyperref}       
\usepackage{url}            
\usepackage{booktabs}       
\usepackage{amsfonts}       
\usepackage{nicefrac}       
\usepackage{microtype}      
\usepackage{lipsum}

\title{A new sigmoidal fractional derivative for regularization}

\author{
   Mostafa Rezapour \\
  Department of Mathematics\\
  Washington State University\\
  Pullman WA, 99163 \\
  \texttt{mostafa.rezapour@wsu.edu} \\
  \And
   Adebowale Sijuwade\\
  Department of Mathematics\\
  Washington State University\\
  Pullman WA, 99163 \\
  \texttt{adebowale.sijuwade@wsu.edu} \\
  \And
   Thomas J. Asaki\\
  Department of Mathematics\\
  Washington State University\\
  Pullman WA, 99163 \\
  \texttt{tasaki@wsu.edu} \\
}

\begin{document}
\maketitle

\begin{abstract}

In this paper, we propose a new fractional derivative, which is based on a Caputo-type derivative with a smooth kernel. We show that the proposed fractional derivative reduces to the classical derivative and has a smoothing effect which is compatible with $\ell_{1}$ regularization. Moreover, it satisfies some classical properties. 

\end{abstract}

\keywords{Fractional calculus, Caputo derivative, Regularization   }

\section{Introduction}
Fractional calculus has undergone significant developments in recent years and has found use in physics, engineering, economics,
 etc \cite{Ref1, Ref2, Ref3}. Classical results about the Riemann-Liouville and Caputo derivatives as well as fractional differential equations can be found in \cite{Ref4, Ref5, Ref6}. In \cite{Ref11} and \cite{Ref48}, Caputo and Fabrizio suggested a new fractional derivative, whose properties were investigated by Losada and Nieto \cite{Ref18}. This fractional derivative was utilized in various applications, including the fractional Nagumo equation in Alqahtani et al. \cite{Ref36},  coupled systems of time-fractional differential problems in Alsaedi et al. \cite{Ref37} and Fischer's reaction-diffusion equation in Atangana et al \cite{Ref38}. More applications of the Caputo-Fabrizio fractional derivative can be found in Aydogan et al \cite{Ref39}. and Atangana et al \cite{Ref40}.

\bigskip 

For $0 \leq \alpha \leq 1$, $ -\infty < a < t$, $f \in H^{1}(a,b)$ and $b>a$, the Caputo fractional derivative is defined by

\[ ^{C} _{a} D_{t}^{\alpha} f(t) =  \frac{1}{\Gamma(1-\alpha)}  \int_{a}^{t} f'(s) (t-s)^{-\alpha} ds. \tag{1} \label{eq: 1} \] 

\noindent By replacing the term $\frac{1}{\Gamma(1-\alpha)}$ with the normalization constant $M(\alpha)$ such that $M(0)=M(1)=1$ and adjusting the kernel $(t-s)^{-\alpha}$, we obtain the Caputo-Fabrizio fractional derivative defined by

\[  ^{CF} _{a} D_{t}^{\alpha} f(t) =  \frac{M(\alpha)}{\Gamma(1-\alpha)} \int_{a}^{t} f'(s) \exp \bigg( \frac{-\alpha(t-s)}{\alpha - 1} \bigg)  ds.   \tag{2} \label{eq: 2}  \] 

\noindent The Caputo-Fabrizio fractional derivative of a constant vanishes as does the usual Caputo derivative, however the new kernel $\exp \big( \frac{-\alpha}{\alpha - 1})$ is no longer singular for $s=t$. Caputo and Fabrizio try to extend their definition in \cite{Ref11} to functions in $L^{1}$ by

\[  ^{CF}  D_{t}^{(\alpha)} f(t) = \frac{M(\alpha)}{\Gamma(1-\alpha)} \int_{-\infty}^{t} (f(s)-f(t)) \exp \bigg( \frac{-\alpha(t-s)}{\alpha - 1} \bigg)  ds.    \]

%
%
%
%
%

Algahtani et al. \cite{Ref36} show that the nonlinear Nagumo equation given by

\[ {^{CF}}_{0} D_{t}^{\alpha} u(x,t) + \beta u(x,t)^{n} \partial_{x} u(x,t) = \partial_{x} (\alpha u(x,t)^{n} \partial_{x} u(x,t)) + \gamma u(x,t) (1-u^{m}) (u^m - \delta),  \tag{3} \label{eq: 3} \]

\noindent where $0<\alpha<1$ and $\beta, \gamma, \delta$ are constant, subject to the boundary conditions 

\[  u(x,0)=f(x), \hspace{1mm}  u(0,t)=g(t) \]

\noindent has an exact solution. The authors show that this PDE can be reformulated in terms of a Lipschitz kernel. Existence of the exact solution is shown using a fixed point approach and uniqueness is provided, given that suitable assumptions are made about the Lipschitz constant. Their study claims that an exponential kernel is in some sense a better kernel than a power function, since the lack of a singularity provides a better filtration effect. In the context of fractional differential equation applications, since the associated functions are not defined in a Banach space, only approximate solutions to certain fractional differential equations can be investigated. The methods used to handle fractional differential problems such as ${^{CF}}D^{\alpha}f(t)= g(t, f(t))$, cannot be extended to the problems resembling ${^{CF}}D^{\alpha}f(t)= g(t, f(t),{^{CF}}D^{\alpha}f(t)$). 

\smallskip

In Baleanu et al \cite{Ref17}, the Caputo-Fabrizio fractional derivative on the Banach space $C_{\mathbb{R}}[0,1]$ is considered in the context of higher order series-type fractional integrodifferential equations. More precisely, an extended Caputo-Fabrizio type fractional derivative is provided of order $0 \leq \alpha < 1$ on  $C_{\mathbb{R}} [0,1]$ for $b>0$ by

\[  {^{CF}} _{N} D^{\alpha} f(t) = \frac{M(\alpha)}{ 1-\alpha} (f(t) - f(0)) \exp \bigg( \frac{-\alpha t}{1-\alpha}   \bigg) + \frac{\alpha M(\alpha)} {   (1-\alpha)^2 } \int_{0}^{t}   (f(t) - f(s))  \exp \bigg( \frac{-\alpha (t-s) }{1-\alpha}   \bigg) ds. \] 

%
%
%
%
%
%
%

\noindent These authors use a standard fixed point approach to establish uniqueness of solutions to fractional series-type differential problems such as

\[  {^{CF}} _{N} D^{\alpha} f(t)  =  \sum\limits_{j=0}^{\infty} \frac{ {^{CF}} _{N} D^{\rho^{[j]}} g(t, f(t), (\phi f)(t), h(t)^{CF} _{N} D^{\gamma}f(t), g(t)^{CF} _{N} D^{\delta} f(t)) 
   }{  2^{j} },  \]  


\noindent with initial condition $f(0)=0$ and  $\alpha, \gamma, \delta, \rho \in (0,1)$.


%


\noindent An extension of this type which is compatible with orders beyond $(0,1)$ has yet to be provided. 

\bigskip

The Caputo-Fabrizio fractional derivative is discussed in the setting of distributions in \cite{Ref41}. Other types of fractional derivatives can be found in Katugampola \cite{Ref35} and Oliveira et al \cite{Ref42}. In de Oliveira \cite{Ref12}, it is shown that the choice of kernel in a Caputo-type fractional derivative is connected to the Laplace transform via convolution.

\bigskip

Let $\mathscr{I}$ denote the Schwarz class of smooth test functions whose derivatives decay at infinity. Moreover, let $\mathscr{I}’$ denote the space of continuous linear functionals on $\mathscr{I}$. The distributional derivative

\noindent \{ T' \} is defined as in \cite{Ref47} 

\[  \int_{\mathbb{R}} T'(t)\phi(t) dt  = -\int_{\mathbb{R}} T(t)\phi(t)dt, \tag{4} \label{eq: 4} \]


\noindent for all smooth compactly supported test functions $\phi$ on $\mathbb{R}$.The distributional Laplace transform is given by 

\[ F(s) = \mathscr{L}(\phi(t))= \mathscr{F} (\phi(t)e^{-\sigma t})(\mu), 
 \]

\noindent where $s = \sigma + i\mu$, $ \mu < 0$ and $\phi(t)e^{-\sigma t} \in \mathscr{I}'$. Suppose that $f$ is supported on $(0, \infty)$ such that $\sigma >0$ and $f(t)e^{-\sigma t} \in \mathscr{I}'$. It follows that the Laplace transform of the derivative is given by

\[ \mathscr{L}(\phi’(t))(s) = s\mathscr{L}(\phi(t)) (s). \] 

\noindent Let $\mathscr{L}$ denote the distributional Laplace transform defined by

\[  \mathscr{L}( f'(x) ) = \mathscr{L}^{-1} ( s \mathscr{L}(f) ). \]

%
%
%
%

\noindent One can define a more general fractional derivative as follows. Suppose that $\Phi(s,\alpha)$ is a fractional integrodifferential operator and $K(t,s): \mathbb{R}^2 \to \mathbb{R}$ is a continuous kernel. Let the corresponding operator $\phi(s,\alpha)$ be defined for some fractional derivative $D^{\alpha}$ such that 

\[ \mathscr{L} ( D^{\alpha} f(t) ) = \Phi(s,\alpha) \mathscr{L}(f(t)), \]


\noindent where $\Phi(s,1)  = s, \Phi(s,-1) = \frac{1}{s}$ and $\Phi(s,0)=1$. Then, letting $\Phi(s, \alpha)=s\mathscr{L}(K(s, t, \alpha))$. Proceeding with the Convolution Theorem, we are left with a Caputo-type fractional operator of the form

\[  _{a} D_{K} ^{\alpha} f(t) = \int_{a}^{t} K(t-s, \alpha) f'(s) ds, \tag{5} \label{eq: 5}  \]

\noindent which is dependent on the choice of kernel $K$. For $f \in H^{1}(a,b)$, and $n \in \mathbb{N}$, we can spot commonly used kernels such as the Caputo kernel $K_{1}= \frac{1}{\Gamma(1-\alpha)} (t-s)^{\ceil{\alpha} - \alpha  - 1}$ , the Caputo-Fabrizio kernel $K_{2} = \frac{M(\alpha)}{1-\alpha} \exp(\frac{-\alpha (t-s)}{\alpha-1})$ and the Gaussian kernel $K_{3}=\frac{1}{\sqrt{2\pi \sigma^2}} \exp(\frac{-t^2}{2\sigma^2})$ \cite{Ref4, Ref10, Ref13}.

\bigskip 

The memory principle for fractional derivatives describes the history of $f(t)$ near the terminal point $t=a$. Let $L$ denote the memory length, satisfying $a+L \leq t \leq b$. Define the error in approximating the fractional derivative by

\[  E_{L,\alpha, a} (t) = |  _{a} D_{K}^{\alpha} f(t)  - _{t-L} D_{K}^{\alpha} f(t) |, \] 

\noindent where $_{a} D_{K} ^{\alpha} f(t)$ is as in ($\ref{eq: 5}$).  If $f'(t) \leq M$ for $a<t<b$ and $0<\alpha<1$, we have the following error estimate for the Caputo fractional derivative

\[  E_{L,\alpha, a} (t) = \bigg|\frac{1}{\Gamma(1-\alpha)} \int_{t-L}^{t} f’(s) (t-s)^{-\alpha} ds\bigg|  \leq \frac{ML^{1-\alpha}}{|\Gamma(2-\alpha)|}.    \] 

\noindent For  all $\epsilon >0$, if $E_{L,\alpha, a} (t) \leq \epsilon$ with $a+L \leq t \leq b$, we have 

\[  L \geq \bigg( \frac{ M }{\epsilon |\Gamma(2-\alpha)|} \bigg)^{\frac{1}{\alpha-1}}. \tag{6} \label{eq: 6} \]

\noindent Therefore, the Caputo fractional derivative with terminal $a$ can be approximated by the corresponding fractional derivative with lower limit $t-L$, with the level of accuracy described above.

\bigskip 

In this work, we propose a different fractional derivative that has a smooth kernel. Our primary interest in defining this fractional derivative is the improvement of machine learning algorithms. Caputo-type fractional derivatives have been applied in machine learning, such as in Pu et al \cite{Ref10}. In particular, fractional order gradient methods have been considered in order to improve the performance of the integer order methods. For example, suppose that $f: \mathbb{R}^{n} \to \mathbb{R}$ is convex and differentiable with a Lipschitz gradient, then the integer order gradient method defined by

\[ x_{k+1} = x_{k}  - \mu \nabla f(x_k) \] 

\noindent has a linear convergence rate. Improving the performance of the integer-order gradient method is critical in optimization problems. In recent literature, fractional calculus has been thought to improve the integer order gradient method due to nonlocality and the memory principle. Fractional order gradient methods have been proposed based on the Caputo fractional derivative that offer competitive convergence rates. For example, in \cite{Ref28}, a Caputo fractional gradient method is proposed that is shown to be monotone and exhibit strong convergence. 

\bigskip

Fractional derivatives were used in the backpropagation algorithm for feedforward neural networks and convolutional neural networks in \cite{Ref32, Ref46}. In both studies, the rate of convergence was shown to exceed the rate of integer-order methods. Fractional-order methods have been used to investigate complex-valued neural networks in \cite{Ref24} and recurrent neural network models in \cite{Ref44}. In \cite{Ref28} and \cite{Ref22}, gradients based on the Caputo fractional derivative are used to update parameters while integer order gradients are used to handle backpropagation allowing for simpler computation. The experiments therein are shown to improve the accuracy of the neural network’s performance compared to integer-order methods while being equally costly.

\bigskip 


In the training of machine learning models, one often needs to obtain weights of the features which optimize the training data. In the case of maximum likelihood training, regularization is typically needed so that the model does not overfit the training data. In $\ell_p$ regularization, the weight vector is penalized by its $\ell_p$ norm. While the case for $p=1$ and $p=2$ are very common and result in similar levels of accuracy, $\ell_1$ regularization is much more practical. Due to its sparsity, $\ell_1$ regularization is less memory intensive and more time-effective than  $\ell_2$ regularization. On the other hand, $\ell_1$ regularization is problematic in that during the update process, the gradient of the regularization term is not differentiable at the origin as the error function given below

\[  E_{\ell_{1}} = E + \lambda \sum\limits_{k=1}^{N} |x_k| \tag{7} \label{eq: 7} \] 

\noindent has classical derivative

\[  \frac{\partial E_{\ell_{1}} } {\partial x_j} = \frac{\partial E }{ \partial x_j }  + \lambda \text{sgn}(x_j). \]

A typical remedy to this problem is to use the stochastic gradient descent method, which approximates the gradient using the training data. Although time efficient for training, when the dimension of the feature space is large, the update process slows down significantly. Furthermore, the model becomes less sparse after training the data. The discontinuity induced by the regularizer proves to be problematic as it adjusts the direction of descent. The use of sigmoids in regularization problems has been previously explored as in Krutikov \cite{Ref43}, but not in the context of fractional derivatives. Another remedy to the aforementioned problem is the use of fractional gradients over the classical descent methods. These methods are still in their infancy and problematic in that convergence to the local optimum is not always guaranteed, even when the algorithm converges as in \cite{Ref9}. Furthermore, these methods often require an adjustment to the fractional derivative by truncation and methods based on memory principle ($\ref{eq: 6}$) due to the computational expense and the failure of the Caputo kernel to be smooth. 

We would also like our operator to be nonlocal. In \cite{Ref13}, it is shown that unlike the Caputo derivative, the Caputo-Fabrizio fractional derivative is not a nonlocal operator. The linear fractional differential equation 

\[ \lambda  (^{CF} _{a} D_{t}^{\alpha} f(t)) + \nu(t)g(t) + \eta(t,t_0)Y(t_0) =0 \] 

\noindent is shown to reduce to a first-order ordinary differential equation. This means that the Caputo-Fabrizio derivative cannot sufficiently describe processes with nonlocality and memory. With the correct choice of kernel, this complication can be avoided.


\bigskip

\section{Main results}

\noindent In this section, we define a new left-sided fractional derivative. We show that the proposed fractional derivative reduces to the $H^{1}$  derivative as the order approaches 1. In the results to follow, for $0<\alpha \leq 1$, we will let $C_1(\alpha)$ denote a normalization constant $\frac{C(\alpha)}{\Gamma(2-\alpha)}$ satisfying $C(\alpha)\Gamma(1-\alpha) \to \frac{1}{2}$ as $\alpha \to 1^{-}$.

\bigskip 

\noindent \textbf{Definition 2.1. (Left sigmoidal fractional derivative)} Let $0<\alpha \leq 1$ , $f \in H^{1}((a,b))$, $t>a$ and \noindent $\{  f(t) \}'$ denotes the $H^{1}$ distributional derivative as in ($\ref{eq: 4}$).  We define a new fractional derivative by

\[  {^{\sigma}} D^{\alpha}_{a} f(t) = C_1(\alpha) \int_{a}^{t} \{ f(s) \}' \sech^2 \bigg( \frac{s-t}{1-\alpha} \bigg) ds. \tag{8} \label{eq: 8} \]

\noindent Now, we show that the left sigmoidal fractional derivative reduces to the $H^{1}$ derivative.

\bigskip 

\begin{theorem}{\textbf{Theorem 2.1. (Reduction to classical derivative) }} Suppose $f \in H^{1}(a,b)$, then

\[  \lim\limits_{\alpha \to 1^{-}} {^{\sigma}} D^{\alpha}_{a} f(t)  = \{  f(t) \}'. \tag{9} \label{eq: 9}  \]

\bigskip \begin{proof} \noindent \textbf{Proof.}  

\[  \lim\limits_{\alpha \to 1^{-}} {^{\sigma}} D^{\alpha}_{a}f(t)   =  \frac{C(\alpha)}{\Gamma(2-\alpha)} \lim\limits_{\alpha \to 1^{-}} \int_{a}^{t} \{ f(s) \}' \sech^2 \bigg( \frac{s-t}{1-\alpha} \bigg) ds \]  

\[ = \frac{2C(\alpha)}{\Gamma(1-\alpha)} \lim\limits_{\alpha \to 1^{-}} \int_{a}^{t} \{ f(s) \}' \frac{ \sech^2 \big( \frac{s-t}{1-\alpha} \big)}{2(1-\alpha)} ds  = \frac{2C(\alpha)}{\Gamma(1-\alpha)} \lim\limits_{\alpha \to 1^{-}} \int_{a}^{t} \{ f(s) \}' \frac{ \sech^2 \big( \frac{s-t}{1-\alpha} \big)}{2(1-\alpha)} ds \]  

\[ =  \lim\limits_{\alpha \to 1^{-}} \frac{2C(\alpha)}{\Gamma(1-\alpha)}  \bigg( \int_{a}^{t} \{ f(s) \}' \lim\limits_{\alpha \to 1^{-}}  \frac{ \sech^2 \bigg( \frac{s-t}{1-\alpha} \bigg)}{2(1-\alpha)} ds \bigg)   = \int_{a}^{t} \{ f(s) \}' \delta(s-t) ds = \{ f(t) \}',         \] 

where the last result follows from the observation that

$\delta(t)$ is the Dirac distribution.

\end{proof}  \end{theorem}

\bigskip 

In the following theorem, we show that this left sigmoidal fractional derivative is commutative with respect to the classical derivative.

\bigskip

\begin{theorem}{\textbf{Theorem 2.2.}} Suppose that $f$ is at least twice continuously differentiable and ${^{\sigma}}  D^{\alpha}_{a} f(t)$ is differentiable. If $f'(a)=0$, then

 \[  {^{\sigma}}  D^{\alpha}_{a} ({^{\sigma}} D^{1}_{a}  f(t)  )   = {^{\sigma}} D^{1}_{a} ({^{\sigma}} D^{\alpha}_{a} f(t)  ) , \tag{10} \label{eq: 10} \]  
where $0<\alpha < 1$.

\bigskip \begin{proof} \noindent \textbf{Proof.}  From ($\ref{eq: 8}$), integrating by parts yields

\[  {^{\sigma}}  D^{\alpha}_{a} ({^{\sigma}} D^{1}_{a}  f(t) )   =  C_1(\alpha) \int_{a}^{t} f''(s) \sech^2 \bigg( \frac{s-t}{1-\alpha} \bigg) ds \]

\[ =   \frac{f'(t)}{1-\alpha}  + \frac{2C(\alpha)}{\Gamma(2-\alpha)(1-\alpha)} \int_{a}^{t} f''(s) \sech( \frac{s-t}{1-\alpha})\tanh( \frac{s-t}{1-\alpha}) ds, \tag{11} \label{eq: 11} \] 

\noindent so we have \[   {^{\sigma}} D^{1}_{a} ({^{\sigma}}  D^{\alpha}_{a} f(t) ) =  \lim\limits_{\gamma \to 1^{-}} {^{\sigma}} D^{\gamma}_{a} ({^{\sigma}}  D^{\alpha}_{a} f(t)  )  = \frac{d}{dt} ({^{\sigma}}  D^{\alpha}_{a} f(t)  )  = C_1(\alpha) \frac{d}{dt}  \int_{a}^{t} f'(s)  \sech^2 \bigg( \frac{s-t}{1-\alpha} \bigg) ds \]

\[ =   \frac{f'(t)}{1-\alpha}  + \frac{2C(\alpha)}{\Gamma(2-\alpha)(1-\alpha)} \int_{a}^{t} f''(s) \sech( \frac{s-t}{1-\alpha})\tanh( \frac{s-t}{1-\alpha}) ds,  \tag{12} \label{eq: 12} \]

appealing to the Leibniz integral rule

\[   \frac{d}{dt} \bigg(  \int_{a(t)}^{b(t)} f(t,s)ds \bigg) = f(t,b(t)) b'(t) - f(t,a(t)) a'(t) + \int_{a(t)}^{b(t)}  \frac{\partial}{\partial t} f(s,t) dt. \]

From ($\ref{eq: 11}$) and ($\ref{eq: 12}$), the desired result is obtained. \end{proof}  \end{theorem}

\bigskip

\noindent  In the next theorem, we show that the left sigmoidal fractional derivative does not satisfy the memory principle in the sense of ($\ref{eq: 6}$). More precisely, the next theorem implies that we show that the left sigmoidal fractional derivative can be approximated by the corresponding fractional derivative with lower limit $t-L$ with increased accuracy for orders in which $C_1(\alpha)$ is large. 
\bigskip

\begin{theorem}{\textbf{Theorem 2.3. (Memory principle)}} Suppose that $f$ is differentiable on $(a,b)$, $ a+L \leq t \leq b$ and $0<\alpha<1$. For every $\epsilon > 0$, if there exists $C_0>0$ such that $f'(t) \leq C_0$, then 

\[ L \geq (1-\alpha)(|C_1(\alpha)|C_0\epsilon^{-1})^{\frac{1}{2}}. \tag{13} \label{eq: 13}    \]

\bigskip \begin{proof} \noindent \textbf{Proof.}  Making use of the inequality 

\[ \cosh(s) \geq \sqrt{1+s^2},\]

\noindent we have

\[  | {^{\sigma}} D^{\alpha}_{a} f(t)  - {^{\sigma}} D^{\alpha}_{t-L} f(t)   |   =  C_1(\alpha) \int_{a}^{t-L}  f'(s) \sech^2 \bigg( \frac{s-t}{1-\alpha} \bigg) ds \leq  C_1(\alpha) C_0 \int_{a}^{t-L} \frac{ds}{1+ ( \frac{s-t}{1-\alpha} )^2 } \leq  \frac{C_1(\alpha)C_0}{ 1 +  ( \frac{L}{1-\alpha} )^2 }, \]

\noindent and the result follows.

\end{proof}  \end{theorem}

In the theorem below, we show that our new fractional derivative provides a sigmoidal approximation to functions that have a piecewise linear $H^{1}$ distributional derivative. For instance, the proposed left sigmoidal fractional derivative is compatible with $\ell_{1}$-regularization. In the case of the $\ell_{1}$ norm, it can be used to define a fractional gradient, which approximates its classical gradient via a family of sigmoids as $\alpha$ approaches 1. This is promising in the context of gradient descent algorithms. 

\bigskip

\begin{theorem}{\textbf{Theorem 2.4 (Norm-1 compatibility})} ${^{\sigma}} D^{\alpha}_{a}$ provides a smooth approximation to the $\ell_{1}$ norm defined by 

\[ \norm{ x }_{1} = \sum\limits_{k=1}^{n} |x_k| \] 

\noindent as $\alpha \to 1$ in the sense that for the error function $E$ given in ($\ref{eq: 7}$), ${^{\sigma}} D^{\alpha}_{a} E_{\ell_{1}}(x_j)$ is given by

\[  {^{\sigma}} D^{\alpha}_{a} E(x_j) + \lambda C_{1}(\alpha)(\alpha - 1) \tanh( \frac{a-x_j}{1-\alpha}), \]  

\noindent where $a>0$.

\bigskip \begin{proof} \noindent \textbf{Proof.}  The result follows from the observation that

\[  C_1(\alpha) \int_{a}^{t} \{ |s| \} \sech^2 \bigg( \frac{s-t}{1-\alpha} \bigg) ds   =  C_1(\alpha) \int_{a}^{t} H(t) \sech^2 \bigg( \frac{s-t}{1-\alpha} \bigg) ds,  \]

\[  =  C_1(\alpha)(\alpha - 1) \tanh( \frac{a-t}{1-\alpha}) ds \to  \frac{1}{2}(2H(t) - 1)  \text{ as } \alpha \to 1^{-}, \]

 \noindent where $H(t)$ is the Heaviside function.

\end{proof}  \end{theorem}

\begin{theorem}{\textbf{Theorem 2.5. (Mittag-Leffler function).} }Suppose that $\gamma, \eta > 0$ and $0<a< t$. Then

\[ {^{\sigma}} D^{\alpha}_{a} E_{\gamma, \eta}(t)  \leq C_1(\alpha) E_{\gamma, \eta}(t-a), \]

\noindent where $E_{\gamma, \eta} (z) =  \sum\limits_{k=0}^{\infty} \frac{z^{k}}{\Gamma(\gamma k + \eta )}$ is the two-parameter Mittag-Leffler function 

\bigskip \begin{proof} \noindent \textbf{Proof.}  \[ {^{\sigma}} D^{\alpha}_{a} E_{\gamma, \eta}(t) =  \int_{a}^{t} \sech^2 \bigg( \frac{s-t}{1-\alpha} \bigg) \frac{d}{ds} \sum\limits_{k=0}^{\infty} \frac{ s^{k}}{\Gamma(\gamma k + \eta)} ds \]

\[ =  \int_{a}^{t} \sech^2 \bigg( \frac{s-t}{1-\alpha} \bigg) \sum\limits_{k=0}^{\infty} \frac{ ks^{k-1}}{\Gamma(\gamma k + \eta)} =    \sum\limits_{k=0}^{\infty} \frac{k }{\Gamma(\gamma k + \eta)} \int_{a}^{t} s^{k-1} \sech^2 \bigg( \frac{s-t}{1-\alpha} \bigg)  ds \]

\[ \leq  \sum\limits_{k=0}^{\infty} \frac{ k }{\Gamma(\gamma k + \eta)} \int_{a}^{t} s^{k-1} ds =  \sum\limits_{k=1}^{\infty} \frac{(t-a)^{k} }{\Gamma(\gamma k + \eta)}. \text{ } \]

\end{proof}  \end{theorem}

\begin{theorem}{\textbf{Theorem 2.6.}} Suppose that $f \geq 0$, $1<p<\infty$, $0<\alpha<1$ and $0< t \leq T$. If $f \geq 0$ is differentiable with $f' \in L^{p}(\mathbb{R})$ and $M$ is the maximal operator of $f$ given by 

\[ Mf(x) = \sup\limits_{t \to 0} \frac{1}{2(a+x)} \int_{a-x}^{a+x} f(t) dt, \]

\noindent then

\begin{enumerate} [label=(\alph*)]

\item ${^{\sigma}} D^{\alpha}_{-t} f(t) \leq 2TC_1(\alpha)M(|f'|)(0)$

\item ${^{\sigma}} D^{\alpha}_{a} f(t)$ is integrable on $\mathbb{R}$.

\end{enumerate}

\bigskip \begin{proof} \noindent \textbf{Proof.}  (a) Since

\[ {^{\sigma}} D^{\alpha}_{-t} Mf(t) = C_1(\alpha) \int_{-t}^{t} f'(s) \sech^2 \bigg( \frac{s-t}{1-\alpha} \bigg) ds \leq  2tC_1(\alpha) \cdot \frac{1}{2t}  \int_{-t}^{t} f'(s) ds \]

\[ \leq 2TC_1(\alpha) \sup\limits_{t>0}  \frac{ \int_{-t}^{t} |f'(s)| ds}{t} = 2TC_1(\alpha) M |f'|(0).  \]

(b) From Young's convolution inequality, $\norm{ f \star g }_{L^{r}} \leq \norm{f}_{L^{p}} \norm{g}_{L^{\frac{pr}{p + r(p-1)}}}$.

\[ \int_{-\infty}^{\infty} {^{\sigma}} D^{\alpha}_{a} Mf(t) dt \leq C_1(\alpha) \int_{-\infty}^{\infty} \int_{a}^{t} |f'(s) \sech^2 \bigg( \frac{s-t}{1-\alpha} \bigg)| ds \] 

\[=\norm{f'(t) \star \sech^2 \bigg( \frac{t}{\alpha-1} \bigg) }_{L^{1}(\mathbb{R})} \leq \norm{ f' }_{L^p (\mathbb{R}) } \norm{\sech^2 \bigg( \frac{t}{\alpha-1} \bigg)}_{L^{\frac{p}{2p-1}}(\mathbb{R})}   < \infty. \text{}\]

\end{proof}  \end{theorem} 

\bigskip 

The next theorem describes the effect of the Laplace and Fourier transforms, which can extend to distributions as in de Oliveira \cite{Ref6}. The Convolution Theorem connects our choice of kernel as in ($\ref{eq: 5}$) via the operator $\Phi(s,\alpha) = s\mathscr{L} (K(s,t,\alpha))$. In this case, $\Phi(s,\alpha)$ depends on the digamma function $\Psi(z)=\frac{\Gamma'(z)}{\Gamma(z)}$. This shows that the left-sigmoidal fractional derivative does not reduce to the left-sided Riemann-Liouville fractional derivative.

\bigskip

\begin{theorem}{\textbf{Theorem 2.7. (Transformations)}} Suppose that $0<\alpha<1$, $Re(s)>0$, $\omega \in \mathbb{R}$, $a \in \mathbb{R}$ and $f$ is a differentiable function of exponential order such that $f(0)=0$. If $T_1(s), T_2(\omega)$ are defined by

\[ T_1(s) = \bigg( 1 + s \bigg( \frac{\Psi(\frac{2 + s}{4}) - \Psi( \frac{s}{4})}{2} \bigg)  \bigg), \hspace{1mm}    T_2(\omega) = \sqrt{ \frac{\pi}{2} } \csch( \frac{\pi \omega}{2} ), \] 

\noindent then

\begin{enumerate} [label=(\alph*)]

\item $ \mathscr{L}( {^{\sigma}} D^{\alpha}_{0} f(t) )(s) = C_1(\alpha)(s(\alpha - 1))^2 {T_1\big({(\alpha - 1)s} \big) \mathscr{L}(f)(s) }$

\item $\mathscr{F}( {^{\sigma}} D^{\alpha}_{0}  )(\omega) = -C_1(\alpha) \omega^2 | {\alpha-1}|(\alpha-1)  {T_2\big( (\alpha - 1)s} \big) \mathscr{F}(f)(\omega), $

\end{enumerate}

\noindent where $\mathscr{L}(f)(s)$ denotes the Laplace transform of $f$ and $\mathscr{F}(f)(\omega)$ denotes the Fourier transform of $f$.

\bigskip \bigskip \begin{proof} \noindent \textbf{Proof.}  (a) follows from a standard application of the Convolution theorem. By using the dilation property $\mathscr{L}(f(at)) = \frac{F(\frac{s}{a})}{a}$, we have

\[ \frac{\mathscr{L} \big( {^{\sigma}} D^{\alpha}_{0} f(t) )(s)}{C_1(\alpha)} = \mathscr{L} \bigg( f' \star \sech^2 \bigg( \frac{t}{\alpha-1} \bigg) \bigg)  = \mathscr{L}(f') \mathscr{L}\bigg(\sech^2 \bigg( \frac{t}{\alpha-1} \bigg) \bigg) \]

\[ = s(\alpha - 1) \mathscr{L}(f)(s)  \mathscr{L}(\sech^2)(s(\alpha-1)) \]

\[ = ((\alpha-1)s)^2 \mathscr{L}(f)(s) \mathscr{L}(\tanh)(s(\alpha-1)).   \]

The transform $\mathscr{L}(\tanh t)$ is handled as follows

\[  s^2\mathscr{L}(\tanh(t))(s)  = s^2 \mathscr{L} ( \tanh (t) )   = s^2 \int_{0}^{\infty} \frac{e^{-st}(1-e^{-2t})}{1 + e^{-2t}} dt  \]

\[ = s^2 \int_{0}^{\infty} e^{-st} (1 - e^{-2t}) \sum\limits_{k=0}^{\infty} (-e^{-2t})^{k}  dt. \]

\noindent Because of the absolute convergence of the monotone decreasing sum $\sum\limits_{k=0}^{\infty} (-1)^ke^{-2kt} dt$ and the nondecreasing nature of its partial sums, we can exchange integration and summation using the Lebesgue Monotone Convergence Theorem. Continuing, we have

\[ s + 2s^2 \sum\limits_{k=1}^{\infty} (-1)^{k} \mathscr{L}(e^{-2kt})  = s + 2s^2 \sum\limits_{k=1}^{\infty} \frac{(-1)^k}{2k + s} \] 

\[ = s + 2s \sum\limits_{k=1}^{\infty}  = \frac{(-1)^k}{\frac{2k}{s} + 1  }  =  s\bigg( 1 + s \bigg( \frac{\Psi(\frac{2 + s}{4}) - \Psi( \frac{s}{4})}{2} \bigg)  \bigg).  \]

\noindent The identity

\[ \sum\limits_{k=0}^{\infty} \frac{(-1)^{k}}{sk + 1} = \frac{ \Psi ( \frac{s+1}{2s}) - \Psi( \frac{1}{2s} )}{2s} \]

\noindent used above comes from the Lerch transcendent, defined by

\[  \Phi(s,z,a) = \sum\limits_{k=0}^{\infty} \frac{z^k}{(a + k )^s}, \]  

\noindent where $|z| < 1$, $a \neq 0, -1, -2 , ... $ and using the dilation property once more, the result follows. 

\bigskip 

\noindent (b) We proceed as in (a). 

\[ \mathscr{F} ( {^{\sigma}} D^{\alpha}_{0} f(t) ) = \int_{-\infty}^{\infty} ( {^{\sigma}} D^{\alpha}_{0} f(t) )  e^{i\omega t} dt = \mathscr{F} (f') \mathscr{F}\bigg( \sech^2\bigg(\frac{t}{\alpha - 1} \bigg) \bigg)  \]

\[ = i \omega \mathscr{F}(f) \mathscr{F}\bigg( \sech^2\bigg(\frac{t}{\alpha - 1} \bigg) \bigg) = i\omega |\alpha - 1| \mathscr{F}(f) \mathscr{F}(\sech^2 )((\alpha - 1)\omega)  \] 

\[  = -\omega^2|\alpha - 1|(\alpha-1) \mathscr{F}(f)  \mathscr{F}(\tanh(t))((\alpha - 1)\omega).  \] 

\noindent To finish the proof, we recall the result 

\[ \mathscr{F}(\tanh (t) )  =  i\omega \sqrt{ \frac{\pi}{2} } \csch( \frac{\pi \omega}{2}) . \text{}\]

\end{proof}  \end{theorem}

\begin{theorem}{\textbf{Theorem 2.8.}} Suppose that $f$ is differentiable and $0<\alpha<1$. Then

\[ \int_{a}^{t} f'(s) e^{ - (\frac{s-t}{1-\alpha})^2}  ds   \leq C_1(\alpha)^{-1} \hspace{0.5mm} {^{\sigma}}  D^{\alpha}_{a} (f(t))  \leq    \int_{a}^{t} \frac{(1-\alpha)^2 f'(s) }{(1-\alpha)^2 + (s-t)^2 }ds \leq f(t) - f(a).  \]  

\bigskip \begin{proof} \noindent \textbf{Proof.}  Using the inequality

\[  \cosh{x} \leq e^{\frac{x^2}{2}}, \]

\noindent we have 

\[   e^{-\frac{1}{2}(\frac{s-t}{1-\alpha})^2} \leq \sech(\frac{s-t}{1-\alpha}), \]

\noindent which results in the leftmost inequality. Noticing that $\cosh^2{x} \geq 1 + x^2$, we have that

\[  \sech^2\bigg( \frac{s-t}{1-\alpha} \bigg) \leq \frac{(1-\alpha)^2 }{(1-\alpha)^2 + (s-t)^2 } \leq 1, \]

\noindent finishing the last three inequalities.

\end{proof}  \end{theorem}

\begin{theorem}{\textbf{Theorem 2.9.}} The problem 

\[ {^{\sigma}}  D^{\alpha}_{a} (f(t))  = G(t), \hspace{1mm} G(0)=0 \] 

\noindent has the solution 

\[  f(t) = \frac{g(t)}{C_1(\alpha)} + f(0), \] 

\noindent where $G(t)=\int_{0}^{t} g(s) ds$.

\bigskip \begin{proof} \noindent \textbf{Proof.}  Differentiating the differential equation above, the problem above reduces to

\[ C_1(\alpha)f'(t) = g'(t),  \]

\noindent which can be integrated to obtain the result.

\end{proof}  \end{theorem}

\begin{theorem}{\textbf{Theorem 2.10.}} Let $0<\alpha<1$ and let $g: (a,b) \times \mathbb{R}^2$ be a continuous function such that there exists a constant $C_0 >0$ satisfying

\[ |g(t,x_1, y_1) - g(t, x_2, y_2)| \leq C_0(|x_1-x_2| + |y_1 - y_2|) \] 

\noindent for all $t \in (a,b)$ and $x_1, x_2, y_1, y_2 \in \mathbb{R}$ and $|(\alpha -1)C(\alpha)C_0| < 1$. Then, the problem

\[  {^{\sigma}}  D^{\alpha}_{a} f(t) = g(t, f(t), {^{\sigma}} D^{\alpha}_{a} f(t) ) \] 

\noindent has a unique solution. 

\bigskip \begin{proof} \noindent \textbf{Proof.}  \[  | g(t, {^{\sigma}}  D^{\alpha}_{a} (f_1(t))) - g(t, {^{\sigma}}  D^{\alpha}_{a} (f_2(t)) |  \] 

\[  \leq  |(\alpha -1)C_1(\alpha) \tanh(\frac{a-t}{1-\alpha}) |f_1 - f_2| \] 

\[  \leq |(\alpha-1)C_1(\alpha) C_0| |f_1 - f_2|. \] 

\noindent Since $(\alpha -1)C(\alpha)C_0 < 1$, the map $F: H^{1}(a,b) \to H^{1} (a,b)$ defined by 

\[ C_1(\alpha)^{-1} g(t, {^{\sigma}}  D^{\alpha}_{a} (f_1(t))) \]

\noindent is a contraction. By the Banach fixed-point theorem, it has a unique fixed point, finishing the proof.

\end{proof}  \end{theorem}

\noindent We note that this result is advantageous in that the analogous existence and uniqueness result as in fractional differential systems defined by the Caputo derivative is highly dependent on initial conditions imposed on the primary function of interest and its classical derivatives\cite{Ref4}.

\bigskip

We now shift our attention to a gradient descent method. Suppose that $f(x)$ has a bounded derivative and unique critical point $t^*$ such that $f’(t^*)=0$. For $a \leq t \leq b$, $0<\alpha<1$, define the scalar left sigmoidal fractional gradient descent method by

\[  t_{k+1} = t_{k} - \mu \hspace{0.5mm} {^\sigma}  D_{t_{k-1}}^{\alpha} f(t_k). \tag{13} \label{eq: 13} \] 

\noindent where $0 < \mu < 1$ is the learning rate.

\bigskip 

\begin{theorem}{\textbf{Theorem 2.11 (Fractional Gradient Descent).}} Let $f$  be as in ($\ref{eq: 13}$). Then, the left-sigmoidal fractional-order gradient method ($\ref{eq: 13}$) converges to the true critical point $t^{*}$. 

\bigskip \begin{proof} \noindent \textbf{Proof.}  Denote the Lipschitz constant of $f$ by $L$. For $k \geq N$, 

\[ |t_{k} - t_{k+1}| = \mu \hspace{0.5mm} {^\sigma}  D_{t_{K-1}}^{\alpha} f(t_k) = C_1(\alpha)\mu \bigg| \int_{t_{k-1}}^{t_k} f’(s) \sech^2 \bigg(  \frac{s-t_k}{1-\alpha}   \bigg) ds \bigg| \]

\[ \leq C_{1}(\alpha) \mu L |\alpha-1| \tanh \bigg( \frac{s-t_k}{1-\alpha} \bigg) \leq C_{1} (\alpha) \mu L |t_{k} - t_{k-1}|. \] 

\noindent Repeating this process, it follows that the $t_k$ form a Cauchy sequence, guaranteeing convergence.  To show that the sequence converges to the critical point, suppose for contradiction that the sequence $(t_k)_{k=0}^{\infty}$ converges to a point $\hat{t} \neq t^{*}$. Then, for every $\epsilon >0$, there exists $N \in \mathbb{N}$ such that for all $k \geq N$, $|f’(t_k)| >0$ and 

\[  | t_{k-1} - \hat{t} | < \epsilon < | t^* - \hat{t} |. \] 

\noindent As a consequence of $($\ref{eq: 13}$)$ and Theorem 2.8, we have

\[ |t_{k+1} - t_{k}| = C_1(\alpha)\mu \bigg| \int_{t_{k-1}}^{t_k} f’(s) \sech^2 \bigg(  \frac{s-a}{1-\alpha}   \bigg) ds \bigg|   \geq C_1(\alpha)\mu \inf\limits_{k > N}  \int_{t_{k-1}}^{t_k}  f'(s) e^{ - (\frac{s-t}{1-\alpha})^2} ds \]

\[ \geq C_1(\alpha)\mu \inf\limits_{k > N} |f’(t_{k-1})|  \int_{t_{k-1}}^{t_k} 1- \bigg( \frac{s-t}{1-\alpha} \bigg)^2 ds \geq M_{1} |t_{k} - t_{k-1}| \bigg( 1 + \frac{ |t_k - t_{k-1}| }{(1-\alpha)^3}  \bigg) \geq   {M_{1} M_{2}|t_{k} - t_{k-1}|^{\frac{3}{2}}}, \]  

\noindent where 

\[  M_{1} = C_1(\alpha)\mu \inf\limits_{k > N} |f’(t_{k-1})| \hspace{1mm}, M_{2} \leq \frac{1}{\sqrt{3(1-\alpha)^3}}. \] 

\noindent On the other hand, we have the inequality

\[  | t_{k+1} - t_{k-1} |  \leq | t_{k+1} - t^{*}| + |t^{*} - t_{k-1}| < 2\epsilon.  \] 

\noindent Choosing $\epsilon <  \frac{1}{2(M_{1} M_{2})^{2}}$ yields $M_1 M_2 > |t_{k+1} - t_{k}|^{\frac{-1}{2}}$, which implies that $|t_{k+1} - t_{k}| > |t_{k} - t_{k-1}|$, contradicting the assumption that the sequence $(t_{k})$ is convergent.

\end{proof}  \end{theorem}

\section{Conclusion}

In this paper, we defined a new sigmoidal fractional derivative, which is compatible with certain weakly differentiable functions. We showed that this fractional derivative satisfies some forms of classical properties and is compatible with the $\ell_{1}$ norm by a sigmoidal approximation. For further research, we will investigate this operator in optimization and machine learning. We note that the left-sigmoidal fractional derivative can be applied in the context of gradient descent, which has applications in optimization and machine learning \cite{Ref7, Ref8}. Recently, backpropagation and convolution neural networks have been studied in the context of fractional derivatives, typically of the Caputo-type are being used for gradient descent. This idea is still novel and needs to see improvements. For example, the gradient descent method has been handled by Sheng et al.; \cite{Ref32}, \cite{Ref33},  Wang et al. ; \cite{Ref28}, Wei et al.; \cite{Ref9} and Bao et al \cite{Ref22}.  These methods are still early in development. The following topics still need to be fully addressed: convergence to an extreme point, extending the available range of fractional order, more complicated neural networks, loss function compatibility and the usage of the chain rule.

\bigskip

\noindent \textbf{Conflict of interest}

\smallskip 

\noindent  The authors declare that there is no conflict of interest.

{}

\end{document}